\documentclass[letterpaper,conference]{ieeeconf}
\usepackage{cite}
\usepackage{color}
\usepackage{graphicx}
\graphicspath{{fig/}{jpeg/}}
\usepackage[cmex10]{amsmath}
\usepackage{amssymb}
\usepackage{algorithm}
\usepackage{algorithmic}
\usepackage{amsmath,amssymb,lipsum}
\usepackage{hyperref}
\usepackage{comment}
\usepackage{graphicx}
\usepackage[font=small]{caption}
\usepackage{subcaption}
\usepackage{setspace}

\input{mysymbol.sty}
\usepackage{needspace}

\newcommand{\myparagraph}[1]{\needspace{1\baselineskip}\medskip\noindent {\it #1.}}

\usepackage{theorem}
\newtheorem{thm}{Theorem}
\newtheorem{lemma}{Lemma}

\newtheorem{corollary}{Corollary}
\newtheorem{assumption}{Assumption}
\newtheorem{remark}{Remark}

\def\eps{\epsilon}

\date{\today}

\IEEEoverridecommandlockouts

\title{\LARGE \bf A Decentralized Second-Order Method for Dynamic Optimization}

\author{Aryan Mokhtari, Wei Shi, Qing Ling, and Alejandro Ribeiro
\thanks{Work supported by NSF CAREER CCF-0952867, ONR N00014-12-1-0997, and NSFC 61004137. A. Mokhtari and A. Ribeiro are with the Dept. of Electrical and Systems Engineering, University of Pennsylvania, Philadelphia, PA 19104, USA. ({\tt aryanm, aribeiro@seas.upenn.edu}). W. Shi is with the Coordinated Science Lab., University of Illinois at Urbana-Champaign, 1308 W Main St, Urbana, IL 61801, USA. ({\tt wilburs@illinois.edu}). Q. Ling is with the Dept. of Automation, University of Science and Technology of China, 96 Jinzhao Rd., Hefei, Anhui, 230026, China. ({\tt qingling@mail.ustc.edu.cn}).} 
}
\begin{document}
\maketitle

%%%%%%%%%%%%%%%%%%%%%%%% MAIN DOCUMENT CONTENT %%%%%%%%%%%%%%%%%%%%%%%%%

\begin{abstract}
This paper considers decentralized dynamic optimization problems where nodes of a network try to minimize a sequence of time-varying objective functions in a real-time scheme. At each time slot, nodes have access to different summands of an instantaneous global objective function and they are allowed to exchange information only with their neighbors. This paper develops the application of the Exact Second-Order Method (ESOM) to solve the dynamic optimization problem in a decentralized manner. The proposed dynamic ESOM algorithm operates by primal descending and dual ascending on a quadratic approximation of an augmented Lagrangian of the instantaneous consensus optimization problem. 
The convergence analysis of dynamic ESOM indicates that a Lyapunov function of the sequence of primal and dual errors converges linearly to an error bound when the local functions are strongly convex and have Lipschitz continuous gradients. Numerical results demonstrate the claim that the sequence of iterates generated by the proposed method is able to track the sequence of optimal arguments.
\end{abstract}


%%%%%%%%%%%%%%%%%%%%%%%%%%%%%%%%
%%%%%%%%%%%%%%%%%%%%%%%%%%%%%%%%
%%%%%     S   E   C   T   I   O   N      %%%%%%%%%%%%%
%%%%%%%%%%%%%%%%%%%%%%%%%%%%%%%%
%%%%%%%%%%%%%%%%%%%%%%%%%%%%%%%%
\section{Introduction}

We consider a decentralized dynamic consensus optimization problem where the components of  a \textit{time-varying} global objective function are available at different nodes of a network. Specifically, consider a discrete time index $t\in\naturals$, a decision variable $\tbx\in \reals^p$, and a connected network containing $n$ nodes where each node $i$ has access to a dynamic local objective $f_{i,t}: \reals^p \to\reals $. The agents' goal is to track the time-varying optimal argument 
%%%
\begin{equation}\label{gen_prob}
\tbx_t^*\ :=\ \argmin_{\tbx\in\reals^p}\ \sum_{i=1}^n f_{i,t}(\tbx),
\end{equation}
%%%
while exchanging information with their neighbors only. Henceforth, we refer to  $f_{i,t}$ as the instantaneous local function of node $i$ at time $t$ and to $\sum_{i=1}^n f_{i,t}$ as the instantaneous aggregate or global objective at time $t$. Distributed dynamic problems like the one in \eqref{gen_prob} are used to formulate problems in distributed signal processing \cite{alriksson2006distributed, farina2010distributed, jakubiec2013d}, distributed control \cite{ogren2004cooperative, borrelli2008distributed,7349151}, and multi-agent robotics \cite{tu2011mobile, zhou2011multirobot, graham2012adaptive}.

For the static version of \eqref{gen_prob} -- with local functions $f_{i,t}=f_i$ that are time invariant and, consequently, with a fixed global objective as well --, there exist numerous descent methods that can solve the problem in a decentralized fashion. Some of these algorithms implement first order descent in the primal domain \cite{Nedic2009,Jakovetic2014-1}, some others rely on first order ascent in the dual domain \cite{rabbat2005generalized, stephanopoulos1975use, jakovetic2015linear, chatzipanagiotis2013augmented, jakovetic2011cooperative}, and some recent efforts attempt to utilize second order information \cite{NN-part1}. Since the dynamic problem in \eqref{gen_prob} can be interpreted as a sequence of static optimization problems, any of the methods in \cite{Nedic2009,Jakovetic2014-1, NN-part1, rabbat2005generalized, stephanopoulos1975use, jakovetic2015linear, chatzipanagiotis2013augmented, jakovetic2011cooperative} can be used as a solution methodology. However, the methods are themselves iterative and their application would require running a large number of (inner) iterations for each of the (outer) time steps $t$; see, e.g., \cite{kar2011gossip}.

Dynamic methods avoid the introduction of multiple time steps and consider that only a few steps of an iterative optimization method are executed for each time index $t$ \cite{braca2010asymptotic, bajovic2011distributed, cattivelli2010diffusion, jakubiec2013d, ling2014decentralized, simonetto2015decentralized, simonettoquasi}. Naturally, these methods track $\tbx_t^*$ with some error because as they implement a descent on $\sum_{i=1}^n f_{i,t}$, the function drifts towards $\sum_{i=1}^n f_{i,t+1}$. These dynamic methods are therefore concerned with characterizing the tracking error \cite{braca2010asymptotic, bajovic2011distributed, cattivelli2010diffusion, jakubiec2013d, ling2014decentralized, simonetto2015decentralized, simonettoquasi} and with developing specific techniques to reduce the steady state gap between the estimated and actual optima \cite{simonetto2015decentralized, simonettoquasi}. Our goal in this paper is to develop the application of the recently proposed exact second order method (ESOM) \cite{mokhtari2016decentralized} for solving the decentralized dynamic optimization problem in \eqref{gen_prob}. 

We begin by introducing decentralized equivalents of \eqref{gen_prob} (Section \ref{sec:problem}) and propose the use of the dynamic ESOM method to solve the resulting decentralized dynamic optimization problem. Dynamic ESOM is a primal-dual algorithm that uses a quadratic approximation of an augmented Lagrangian (Section \ref{sec:ESOM}). This approximation is expected to have good convergence properties because it incorporates second order information. Alas, this quadratic approximation requires access to the Hessian inverse of the augmented Lagrangian, which is not locally computable. This issue is resolved by using a truncation of the Taylor's series expansion of the Hessian inverse 
%\cite{ZarghamEtal14,NN-part1,NN-part2}
\cite{NN-part1} (Section \ref{sec:hessian_approx}). We study convergence properties of dynamic ESOM and show that the sequence of iterates it generates converges linearly to a neighborhood of the sequence of optimal arguments $\bbx_t^*$ (Section \ref{sec_convg}). We perform a numerical evaluation of the performance of dynamic ESOM in solving a dynamic least squares problem (Section \ref{sec:num_exp}) and close the paper with concluding remarks (Section \ref{sec:conclusions}). 

\myparagraph{\bf Notation} Vectors are written as
$\bbx\in\reals^p$ and matrices as $\bbA\in\reals^{p\times p}$.
Given $n$ vectors $\bbx_i$, the vector
$\bbx=[\bbx_1;\ldots;\bbx_n]$ represents a stacking of the
elements of each individual $\bbx_i$. We use $\|\bbx\|$ and  $\|\bbA\|$ to denote
the Euclidean norm of vector $\bbx$ and matrix $\bbA$, respectively. The norm of vector $\bbx$ with respect to positive definite matrix $\bbA$ is $\|\bbx\|_\bbA:=(\bbx^T\bbA\bbx)^{1/2}$. Given a function $f$ its gradient evaluated at $\bbx$ is denoted as $\nabla f(\bbx)$ and its Hessian as $\nabla^2 f(\bbx)$. The diagonalized version of matrix $\bbA$ is denoted by $\diag(\bbA)$ where its diagonal components are identical with those of $\bbA$ and the other components are null.

%%%%%%%%%%%%%%%%%%%%%%%%%%%%%%%%
%%%%%%%%%%%%%%%%%%%%%%%%%%%%%%%%
%%%%%     S   E   C   T   I   O   N      %%%%%%%%%%%%%
%%%%%%%%%%%%%%%%%%%%%%%%%%%%%%%%
%%%%%%%%%%%%%%%%%%%%%%%%%%%%%%%%
\section{Problem formulation}\label{sec:problem}
Consider $\bbx_i\in\reals^p$ as the copy of the decision variable $\tbx$ at node $i$ and define $\ccalN_{i}$ as the neighborhood of node $i$. Connectivity of the network implies that problem \eqref{gen_prob} is equivalent to the optimization problem
%%%
\begin{align}\label{original_optimization_problem2}
   \{\bbx_{i,t}^*\}_{i=1}^n\ := \
   &\argmin_{\{\bbx_{i}\}_{i=1}^n} \ \sum_{i=1}^{n}\ f_{i,t}(\bbx_{i}), \nonumber\\ 
   &\text{\ s.t.}  \ \bbx_{i}=\bbx_{j}, 
                   \quad \text{for all\ } i, j\in\ccalN_i .
\end{align} 
%%%%
To verify the equivalence of \eqref{gen_prob} and \eqref{original_optimization_problem2}, note that a set of feasible solutions for \eqref{original_optimization_problem2} has the general form of $\bbx_1=\dots=\bbx_n$, since the network is connected. Likewise, the optimal solution of \eqref{original_optimization_problem2} satisfies $\bbx_{1,t}^*=\dots=\bbx_{n,t}^*$. When the arguments $\bbx_i$ of the functions $f_{i,t}(\bbx_{i})$ are equal to each other the objective function $\sum_{i=1}^{n} f_{i,t}(\bbx_{i})$ in \eqref{original_optimization_problem2} can be simplified as the aggregate function $\sum_{i=1}^{n} f_{i,t}(\bbx)$ in \eqref{gen_prob}. Thus, the optimal argument of each node $\bbx_{i,t}^*$ in \eqref{original_optimization_problem2} is identical to the optimal solution $\tbx_t^*$ of \eqref{gen_prob}, i.e., $\bbx_{1,t}^*=\dots=\bbx_{n,t}^*=\tbx_t^*$.

To derive the update for the dynamic ESOM algorithm, define $\bbx:=[\bbx_1;\dots;\bbx_n]\in \reals^{np}$ as the concatenation of the local decision variables $\bbx_i$ and the global function $f_t:\reals^{np}\to\reals$ at time $t$ as $f_t(\bbx)=f_t(\bbx_1,\dots,\bbx_n):=\sum_{i=1}^n f_{i,t}(\bbx_i)$. Further, we introduce the weight matrix $\bbW\in\reals^{n\times n}$ where the element $w_{ij}\geq 0$ represents the weight that node $i$ assigns to node $j$. The weight $w_{ij}$ is nonzero if and only if $j\in \ccalN_i$ or $j=i$. We assume that the assigned weights are chosen such that the weight matrix $\bbW$ satisfies the following conditions
 \begin{equation}\label{weight_ass}
 \bbW=\bbW^T, \quad \bbW\bbone=\bbone, \quad \text{null}(\bbI-\bbW)= \text{span}(\bbone).
 \end{equation}
%%%
The first condition $ \bbW=\bbW^T$ implies that the weights are symmetric, i.e., $w_{ij}=w_{ji}$. The condition $\bbW\bbone=\bbone$ ensures that the weight matrix $\bbW$ is doubly stochastic and the matrix $\bbI-\bbW$ has a zero eigenvalue where its corresponding eigenvector is vector $\bbone$. The last condition $\text{null}(\bbI-\bbW)= \text{span}(\bbone)$ ensures that the matrix $\bbI-\bbW$ has rank $n-1$ and the condition $(\bbI-\bbW)\bbv=\bb0$ holds if and only if $\bbv\in \text{span}\{\bbone\}$. Conditions in \eqref{weight_ass} are typical of mixing matrices and they are required to enforce consensus. 
  
It has been shown (Proposition 1 in \cite{mokhtari2016decentralized}), if we define the matrix $\bbZ=\bbW\otimes \bbI_p\in \reals^{np\times np}$ as the Kronecker product of the weight matrix $\bbW$ and the identity matrix $\bbI_p$, the optimization problem in \eqref{original_optimization_problem2} can be written as
\begin{equation}\label{constrained_opt_problem}
\bbx_t^* =\argmin_{\bbx\in \reals^{np}} \ f_t(\bbx)\qquad \text{s.t.}\ (\bbI-\bbZ)^{1/2}\bbx=\bb0.
\end{equation}
Thus, the optimization problem in \eqref{constrained_opt_problem} is equivalent to the original dynamic problem in \eqref{gen_prob} and we proceed to develop dynamic ESOM to solve \eqref{constrained_opt_problem} in lieu of \eqref{gen_prob}. By introducing $\bbv\in\reals^{np}$ as the dual variable associated with the constraint $ (\bbI-\bbZ)^{1/2}\bbx=\bb0$ in \eqref{constrained_opt_problem}, we define the augmented Lagrangian $\ccalL_t(\bbx,\bbv)$ of \eqref{constrained_opt_problem} as
\begin{equation}\label{lagrangian}
\ccalL_t(\bbx,\bbv)= f_t(\bbx)+\bbv^T(\bbI-\bbZ)^{1/2}\bbx+\frac{\alpha}{2}\bbx^T(\bbI-\bbZ)\bbx,
\end{equation}
where $\alpha$ is a positive constant. Based on the properties of the matrix $\bbZ$, the inner product $\bbx^T(\bbI-\bbZ)\bbx$ augmented to the Lagrangian is null when the variable $\bbx$ is a feasible solution of \eqref{constrained_opt_problem}, otherwise the inner product is positive and behaves as a penalty for the violation of the consensus constraint. 

A well studied approach to estimate the instantaneous minimizer $\bbx_t^*$ is to define $\bbx_t$ as the minimizer of the proximal augmented Lagrangian which is the sum of the augmented Lagrangian $\ccalL_t(\bbx,\bbv_{t-1})$ and the proximal term $(\eps/2)\|\bbx-\bbx_{t-1}\|^2$. This scheme can be interpreted as a dynamic extension of the proximal method of multipliers \cite{hestenes1969multiplier,bertsekas2014constrained}.  Thus, the estimator $\bbx_t$ is the minimizer of the optimization problem
\begin{equation}\label{PMM_primal_update}
\bbx_{t}=\argmin_{\bbx \in \reals^{np}} \left\{\ccalL_t(\bbx,\bbv_{t-1})+\frac{\eps}{2}\|\bbx-\bbx_{t-1}\|^2\right\},
\end{equation}
where $\bbv_{t-1}$ is the dual variable evaluated at step $t-1$ and $\eps$ is a positive constant. The updated dual variable $\bbv_t$ is updated by ascending through the augmented Lagrangian gradient $\nabla_{\bbv}\ccalL_t(\bbx_t,\bbv_{t-1})$ with respect to $\bbv$ with stepsize $\alpha$,
\begin{equation}\label{PMM_dual_update}
\bbv_{t}=\bbv_{t-1} + \alpha (\bbI-\bbZ)^{1/2}\bbx_{t}.
\end{equation}

However, there are two issues with the updates in \eqref{PMM_primal_update}. The first issue is the computation time of the update, since the minimization could be computationally costly. The second drawback is the quadratic term $\bbx^T(\bbI-\bbZ)\bbx$ in \eqref{PMM_primal_update} which is not separable. Thus, the update is not implementable in a decentralized fashion. To resolve these issues we introduce the dynamic ESOM algorithm in the following section.

%%%%%%%%%%%%%%%%%%%%%%%%%%%%%%%%
%%%%%%%%%%%%%%%%%%%%%%%%%%%%%%%%
%%%%%       S   E   C   T   I   O   N       %%%%%%%%%%%%
%%%%%%%%%%%%%%%%%%%%%%%%%%%%%%%%
%%%%%%%%%%%%%%%%%%%%%%%%%%%%%%%%
\section{Dynamic ESOM}\label{sec:ESOM}

In this section, we introduce the dynamic ESOM algorithm as a decentralized algorithm that replaces the augmented Lagrangian $\ccalL_t(\bbx,\bbv_{t-1})$ in \eqref{PMM_primal_update} by its quadratic approximation. This modification reduces the computational complexity of the update in \eqref{PMM_primal_update} and leads to a separable primal update. In particular, we approximate the augmented Lagrangian $\ccalL_t(\bbx,\bbv_{t-1})$ in \eqref{PMM_primal_update} by its second-order Taylor's expansion near the point $(\bbx_{t-1},\bbv_{t-1})$ which is given by
$\ccalL_t(\bbx_{t-1},\bbv_{t-1}) + \nabla_{\bbx}\ccalL_t(\bbx_{t-1},\bbv_{t-1})^T(\bbx-\bbx_{t-1}) + ({1/2})(\bbx-\bbx_{t-1})^T\nabla_{\bbx\bbx}^2\ccalL_t(\bbx_{t-1},\bbv_{t-1})(\bbx-\bbx_{t-1})$. Applying this substitution leads to the update
\begin{align}\label{ESOM_primal_update}
&\bbx_{t} = \\
 & \argmin_{\bbx \in \reals^{np}}  \Big\{\ccalL_t(\bbx_{t-1},\bbv_{t-1})+ \nabla_{\bbx}\ccalL_t(\bbx_{t-1},\bbv_{t-1})^T\!(\bbx-\bbx_{t-1}) \nonumber \\
&\quad + \frac{1}{2}(\bbx-\bbx_{t-1})^T\!\!\left(\nabla_{\bbx\bbx}^2\ccalL_t(\bbx_{t-1},\bbv_{t-1})+\eps \bbI\right)\!(\bbx-\bbx_{t-1})\Big\}.\nonumber
\end{align}
Solving the minimization in the right hand side of \eqref{ESOM_primal_update}
and using the definition of the augmented Lagrangian $\ccalL_t(\bbx,\bbv)$ in \eqref{lagrangian}, it follows that the variable $\bbx_{t}$ can be evaluated as
%%%
\begin{align}\label{ESOM_primal_update_2}
\bbx_{t}=\bbx_{t-1} - \bbH_{t}^{-1}\bigg[\nabla f_t(\bbx_{t-1})&+ (\bbI-\bbZ)^{1/2}\bbv_{t-1}
\nonumber\\
&\quad
 +{\alpha}(\bbI-\bbZ)\bbx_{t-1}\bigg],
\end{align}
where the matrix $\bbH_{t}\in \reals^{np\times np}$ is defined as the Hessian of the objective function in \eqref{ESOM_primal_update} which is given by
\begin{equation}\label{exact_Hessian}
\bbH_{t}:=\nabla^2f_t(\bbx_{t-1})+{\alpha}(\bbI-\bbZ)+\eps\bbI.
\end{equation}
The Hessian $\bbH_{t}$ in \eqref{exact_Hessian} is a block neighbor sparse matrix. In other words, its $(i,j)$th block, which is in $\reals^{p\times p}$, is non-zero if and only if $j\in\ccalN_i$ or $j=i$. This is true since the matrix $\nabla^2f_t(\bbx_{t-1})+\eps \bbI$ is block diagonal and the matrix ${\alpha}(\bbI-\bbZ)$ is block neighbor sparse. Albeit, the Hessian $\bbH_{t}$ is block neighbor sparse, its inverse $\bbH_{t}^{-1} $ in \eqref{ESOM_primal_update_2} is not. Thus, the nodes cannot implement the update in \eqref{ESOM_primal_update_2} in a decentralized fashion.

 To resolve this issue, we use a Hessian inverse approximation that is built on truncating the Taylor's series of the Hessian inverse $\bbH_{t}^{-1} $ as in \cite{NN-part1}. To be precise, we decompose the Hessian as $\bbH_{t}=\bbD_{t}-\bbB$ where $\bbD_{t}$ is a block diagonal positive definite matrix and $\bbB$ is a neighbor sparse positive semidefinite matrix. We define the matrix $\bbD_{t}$ as
\begin{equation}
\bbD_{t} := \nabla^2f_t(\bbx_{t-1})+\eps \bbI+ 2{\alpha}(\bbI-\bbZ_d),
\end{equation}
%%%
where $\bbZ_d:=\diag(\bbZ)$. Hence, the relation $\bbB=\bbD_{t}-\bbH_{t}$ implies that 
\begin{equation}
\bbB :={\alpha}\left(\bbI-2\bbZ_d+\bbZ\right).
\end{equation}
%%%
Considering the decomposition $\bbH_{t}=\bbD_{t}-\bbB$, it follows that the Hessian inverse $\bbH_{t}^{-1}=(\bbD_{t}-\bbB)^{-1}$ can be written as $\bbH_{t}^{-1}=\bbD_{t}^{-1/2}(\bbI-\bbD_{t}^{-1/2}\bbB\bbD_{t}^{-1/2})^{-1}\bbD_{t}^{-1/2}$ by factoring $\bbD_{t}^{1/2}$ from both sides. Note that the absolute value of the eigenvalues of the matrix $\bbD_{t}^{-1/2}\bbB\bbD_{t}^{-1/2}$ are strictly smaller than $1$; see e.g. Proposition 2 in \cite{NN-part1}. Thus, we can use the Taylor's series $(\bbI-\bbX)^{-1}=\sum_{u=0}^{\infty} \bbX^u$ for $\bbX=\bbD_{t}^{-1/2}\bbB\bbD_{t}^{-1/2}$ to write the Hessian inverse $\bbH_{t}^{{-1}}$ as
%%%
\begin{equation}\label{Hessian_approx00}
\bbH_{t}^{{-1}}:=\bbD_{t}^{-1/2}\
\sum_{u=0}^\infty\left(\bbD_{t}^{-1/2}\bbB\bbD_{t}^{-1/2}\right)^u\
\bbD_{t}^{-1/2}.
\end{equation}
%%%
Computation of the Hessian inverse $\bbH_{t}^{{-1}}$ in \eqref{Hessian_approx00} requires infinite rounds of communication between the nodes; however, we can approximate the Hessian inverse $\bbH_{t}^{-1}$ by truncating the first $K+1$ terms of the sum in \eqref{Hessian_approx00}. This approximation leads to the Hessian inverse approximation
%%%
\begin{equation}\label{Hessian_approx}
\hbH_{t}^{{-1}}(K):=\bbD_{t}^{-1/2}\
\sum_{u=0}^K\left(\bbD_{t}^{-1/2}\bbB\bbD_{t}^{-1/2}\right)^u\
\bbD_{t}^{-1/2}.
\end{equation}
%%%
The approximate Hessian inverse $\hbH_{t}^{{-1}}(K)$ is $K$-hop block neighbor sparse, i.e., its $(i,j)$th block is nonzero if and only if there exists at least one path between nodes $i$ and $j$ of length $K$ or smaller. 

We introduce the dynamic ESOM algorithm as a second-order method for solving the decentralized consensus optimization problem which substitutes the Hessian inverse $\bbH_{t}^{{-1}}$ in \eqref{ESOM_primal_update_2} by the $K$-hop block neighbor sparse Hessian inverse approximation $\hbH_k^{-1}(K)$ defined in \eqref{Hessian_approx}. Thus, the update for the primal variable of dynamic ESOM is given by
%%%
\begin{align}\label{ESOM_primal_update_3}
\bbx_{t}=\bbx_{t-1} - \hbH_{t}^{{-1}}(K)&\bigg[\nabla f_t(\bbx_{t-1})+(\bbI-\bbZ)^{1/2}\bbv_{t-1}
\nonumber\\
&\qquad\qquad \
+{\alpha}(\bbI-\bbZ)\bbx_{t-1}\bigg].
\end{align}
The update for the dual variable $\bbv_t$ of dynamic ESOM is identical to the update in \eqref{PMM_dual_update},
%%%
%%
\begin{equation}\label{ESOM_dual_update}
\bbv_{t}=\bbv_{t-1} +\alpha (\bbI-\bbZ)^{1/2}\bbx_{t}.
\end{equation}

Note that the primal and dual updates of dynamic ESOM in \eqref{ESOM_primal_update_3} and \eqref{ESOM_dual_update} are different from the updates of ESOM in \cite{mokhtari2016decentralized} which is designed for static consensus optimization. In particular, the primal and dual updates of ESOM are derived by approximating the time-invariant augmented Lagrangian $\ccalL(\bbx,\bbv):=f(\bbx)+\bbv^T(\bbI-\bbZ)^{1/2}\bbx+({\alpha/2})\bbx^T(\bbI-\bbZ)\bbx$, while the updates for dynamic ESOM are established by a quadratic approximation of the time-variant augmented Lagrangian $\ccalL_t(\bbx,\bbv)$ defined in \eqref{lagrangian}.

The updates in \eqref{ESOM_primal_update_3} and \eqref{ESOM_dual_update} explain the rationale behind dynamic ESOM; however, they are not implementable in a decentralized fashion, since the squared matrix $(\bbI-\bbZ)^{1/2}$ is not block neighbor sparse. In the following section, we introduce a new set of updates for dynamic ESOM which are implementable in a distributed fashion, while they are equivalent to the updates in \eqref{ESOM_primal_update_3} and \eqref{ESOM_dual_update}.

%%%%%%%%%%%%%%%%%%%%%%%%%%%%%%%%
%%%%%%%%%%%%%%%%%%%%%%%%%%%%%%%%
%%%%%    S   U   B  --  S   E   C   T   I   O   N      %%%%%%%
%%%%%%%%%%%%%%%%%%%%%%%%%%%%%%%%
%%%%%%%%%%%%%%%%%%%%%%%%%%%%%%%%
\subsection{Decentralized implementation of dynamic ESOM}\label{sec:hessian_approx}

To come up with updates for dynamic ESOM that can be implemented in a decentralized setting, define the sequence of variables $\bbq_t$ as $\bbq_{t}:=(\bbI-\bbZ)^{1/2}\bbv_t$. Substitute the term $(\bbI-\bbZ)^{1/2}\bbv_t$ in \eqref{ESOM_primal_update_3} by $\bbq_t$ to rewrite the primal update as
%%%
\begin{align}\label{ESOM_primal_update_4}
\bbx_{t}=\bbx_{t-1} - \hbH_{t}^{{-1}}(K)\!\left[\nabla f_t(\bbx_{t-1})+ \bbq_{t-1}+{\alpha}(\bbI-\bbZ)\bbx_{t-1}\right]\!.
\end{align}
Multiplying the dual update in \eqref{ESOM_dual_update} by $(\bbI-\bbZ)^{1/2}$ from the left hand side and using the definition $\bbq_{t}:=(\bbI-\bbZ)^{1/2}\bbv_t$, it follows that
\begin{equation}\label{ESOM_dual_update_4}
\bbq_{t}=\bbq_{t-1} +\alpha(\bbI-\bbZ)\bbx_{t}.
\end{equation}
The system of updates in \eqref{ESOM_primal_update_4} and \eqref{ESOM_dual_update_4} are implementable in a decentralized fashion, since the matrix $\bbI-\bbZ$, which is required for both updates, is block neighbors sparse. Notice that the updates in \eqref{ESOM_primal_update_4} and \eqref{ESOM_dual_update_4} are equivalent to the updates in  \eqref{ESOM_primal_update_3} and \eqref{ESOM_dual_update}, i.e., the sequence of iterates $\bbx_t$ generated by these two schemes are identical. 

We proceed to derive the local updates at each node to implement the primal and dual updates in \eqref{ESOM_primal_update_4} and \eqref{ESOM_dual_update_4}, respectively. To do so, define $\bbg_{t}$ as the augmented Lagrangian gradient $\nabla_\bbx \ccalL_t(\bbx_{t-1},\bbv_{t-1})$ with respect to $\bbx$ which is given by
%%%
\begin{equation}\label{gradient}
\bbg_{t}=\nabla f_t(\bbx_{t-1})+ \bbq_{t-1}+{\alpha}(\bbI-\bbZ)\bbx_{t-1}
\end{equation}
%%%
Further, define the primal descent direction $\bbd_{t}^{(K)}$ evaluated using the Hessian inverse approximation $\hbH_{t}^{{-1}}(K)$ with $K$ levels of approximation as 
\begin{equation}\label{des_new}
\bbd_{t}^{(K)}:=-\hbH_{t}^{{-1}}(K) \bbg_{t}.
\end{equation}
The definition of the descent direction $\bbd_{t}^{(K)}$ in \eqref{des_new} allows us to rewrite the update in \eqref{ESOM_primal_update_4} as $\bbx_{t}=\bbx_{t-1}+\bbd_{t}^{(K)}$. According to the mechanism of Hessian inverse approximation in \eqref{Hessian_approx}, the descent directions $\bbd_{t}^{(k)}$ and $\bbd_{t}^{(k+1)}$ satisfy the recursion
%%%
\begin{equation}\label{relation}
\bbd_{t}^{(k+1)}= \bbD_{t}^{-1}\bbB \bbd_{t}^{(k)}-\bbD_{t}^{-1}\bbg_{t}.
\end{equation}
%%%
Consider $\bbd_{i,t-1}^{(k)}$ as the descent direction of node $i$ at step $t$ which is the $i$-th element of the global descent direction $\bbd_{t}^{(k)}=[\bbd_{1,t}^{(k)};\dots;\bbd_{n,t}^{(k)}]$. Use this definition to write the localized version of the relation in \eqref{relation} at node $i$ as
%%%
\begin{equation}\label{relation_local}
\bbd_{i,t}^{(k+1)}= \bbD_{ii,t}^{-1}\sum_{j=i,j\in\ccalN_i}\left(\bbB_{ij} \bbd_{j,t}^{(k)}\right)-\bbD_{ii,t}^{-1}\bbg_{i,t},
\end{equation}
%%%
where $\bbD_{ii,t}$ is the $i$-th diagonal block of the matrix $\bbD_t$ and $\bbB_{ij}$ is the $(i,j)$-th block of the matrix $\bbB$. 
 Based on the expression in \eqref{relation_local}, node $i$ is able to compute its descent direction $\bbd_{i,t}^{(k+1)}$ using the $k$-th level descent direction of itself $\bbd_{i,t}^{(k)}$ and its neighbors $\bbd_{j,t}^{(k)}$ for $j\in \ccalN_i$. Therefore, if nodes exchange their $k$-th level descent direction $\bbd_{i,t}^{(k)}$ with their neighbors, they can compute the $(k+1)$-th level descent direction $\bbd_{i,t}^{(k+1)}$.

Notice that the block $\bbD_{ii,t}:=\nabla^2 f_{i,t}(\bbx_{i,t-1}) + (2\alpha(1-w_{ii}))\bbI+\eps\bbI$ is locally available at node $i$. Moreover, node $i$ can evaluate the blocks $\bbB_{ii}=\alpha(1-w_{ii})\bbI$ and $\bbB_{ij}=\alpha w_{ij}\bbI$ locally. In addition, nodes can compute the gradient $\bbg_t$ by communicating with their neighbors. To confirm this claim, observe that the $i$-th element of the gradient $\bbg_t=[\bbg_{1,t};\dots;\bbg_{n,t}]$ associated with node $i$ is given by
%%%
\begin{align}\label{local_gradient}
\bbg_{i,t}&:=\nabla f_t(\bbx_{i,t-1})+\bbq_{i,t-1}+\alpha(1-w_{ii})\bbx_{i,t-1}\nonumber\\
&\qquad-\alpha\sum_{j\in \ccalN_i}{w_{ij}}\bbx_{j,t-1},
\end{align}
%%%
where $\bbq_{i,t-1}\in \reals^{p}$ is the $i$-th element of the vector $\bbq_{t-1}=[\bbq_{1,t-1};\dots;\bbq_{n,t-1}]\in \reals^{np}$. Hence, node $i$ can compute its local gradient $\bbg_{i,t}$ using local information $\bbq_{i,t-1}$ and $\bbx_{i,t-1}$, and its neighbors' information $\bbx_{j,t-1}$ where $j\in\ccalN_i$. 

The recursive update in \eqref{relation_local} shows that at each step $t$ nodes can compute the descent directions $\bbd_{i,t}^{(K)}$ by $K$ rounds of communication with their neighbors. Moreover, observe that nodes can implement the dual update in \eqref{ESOM_dual_update_4} as 
\begin{equation}\label{dual_local}
\bbq_{i,t}=\bbq_{i,t-1} +\alpha(1-w_{ii})\bbx_{i,t}-\alpha \sum_{j\in\ccalN_i} w_{ij}\bbx_{j,t},
\end{equation}
by having access to the updated primal variables $\bbx_{j,t}$ of their neighbors $j\in\ccalN_i$.

%%%%%%%%%%%%%%%%%%%%%%%%%%%%%%%%%%%%
%%%%%%%%%%%%%%%%%%%%%%%%%%%%%%%%%%%%
%%%   A   L   G   O   R   I   T   H   M    %%%%%%%%%%%%%%%%%
%%%%%%%%%%%%%%%%%%%%%%%%%%%%%%%%%%%%
%%%%%%%%%%%%%%%%%%%%%%%%%%%%%%%%%%%%
\begin{algorithm}[t]{\small
\caption{Dynamic ESOM-$K$ method at node $i$}\label{algo_ESOM} 
\begin{algorithmic}[1] {
\REQUIRE  Initial iterates $\bbx_{i,0}\!=\!\bbx_{j,0}\!=\!\bb0$ for $j\!\in\!\ccalN_i$ and $\bbq_{i,0}\!=\!\bb0$.
\STATE Compute $\bbB_{ii}=\alpha(1-w_{ii})\bbI_p$ and $\bbB_{ij}=\alpha w_{ij}\bbI_p$ all $j\in\ccalN_i$
\FOR {times $t=1,2,\ldots$}
\STATE Observe the local function $f_{i,t}$
   \STATE Compute $\bbD_{ii,t}= \nabla^2 f_{i,t}(\bbx_{i,t-1}) +2\alpha(1-w_{ii})\bbI_p $
   \STATE Compute the gradient   
         $ \bbg_{i,t} $ as in \eqref{local_gradient}
          %inner loop
   \STATE Compute the initial descent direction $\bbd_{i,t}^{(0)}=-\bbD_{ii,t}^{-1}\bbg_{i,t}$\\ 
   \FOR  {$k=  0, \ldots, K-1$ } 
      \STATE Exchange $\bbd_{i,t-1}^{(k)}$ with neighbors $j\in\ccalN_i$
      \STATE Compute $\displaystyle{\bbd_{i,t}^{(k+1)}= \bbD_{ii,t}^{-1}\!\!\sum_{j=i,j\in\ccalN_i}\!\!\bbB_{ij} \bbd_{j,t}^{(k)}-\bbD_{ii,t}^{-1}\bbg_{i,t}}$
   \ENDFOR
  %  $\bbD_{t}_{N}(t)=$ NN Descent Direction $\left(\bbD_{t}_{0}(t), \bbg(t), N \right)$.
          %
          \STATE Update primal iterate: 
          $\displaystyle{\bbx_{i,t}=\bbx_{i,t-1} +\ \bbd_{i,t-1}^{(k)}}$.
           \STATE Exchange iterates $\bbx_{i,t}$ with neighbors $\displaystyle{j\in \mathcal{N}_i}$.
           \STATE Update the dual iterate:\\ $\displaystyle{\bbq_{i,t}=\bbq_{i,t-1} +\alpha(1-w_{ii})\bbx_{i,t}-\alpha \sum_{j\in\ccalN_i} w_{ij}\bbx_{j,t}}$
\ENDFOR}
\end{algorithmic}}\end{algorithm}

%%%%%%%%%%%%%%%%%%%%%%%%%%%%%%%%
%%%%%%%%%%%%%%%%%%%%%%%%%%%%%%%%
%%%         M  A  I  N        M  A  T  T  E  R     %%%%%%%%%%
%%%%%%%%%%%%%%%%%%%%%%%%%%%%%%%%
%%%%%%%%%%%%%%%%%%%%%%%%%%%%%%%% 
The steps of dynamic ESOM-$K$ at node $i$ are summarized in Algorithm \ref{algo_ESOM}. In step 3, each node $i$ observes its local function $f_{i,t}$ for the current time $t$ and uses this information to compute the block $\bbD_{ii,t}$ and the local gradient $\bbg_{i,t}$ in Steps 4 and 5, respectively. Node $i$ computes its $(k+1)$-th level descent direction $\bbd_{i,t}^{(k+1)}$ in Step 9 using the $k$-th level local descent direction $\bbd_{i,t}^{(k)}$ and the neighbors' decent directions $\bbd_{j,t}^{(k)}$ which are exchanged in Step 8. Note that the recursion in Steps 8 and 9 are initialized by the descent direction $\bbd_{i,t}^{(0)}$ of dynamic ESOM-$0$ evaluated in Step 6. Each node computes its local primal variable $\bbx_{i,t}$ in Step 11 and exchanges it with its neighbors in Step 12. The dual variables $\bbq_{i,t}$ can be updated in Step 13, using the updated local and neighboring  primal variables. The blocks $\bbB_{ij}$ for $j=\ccalN_i$ and $j=i$ are time invariant and they are computed and stored locally in Step 1. 

%%%%%%%%%%%%%%%%%%%%%%%%%%%%%%%%
%%%%%%%%%%%%%%%%%%%%%%%%%%%%%%%%
%%%         R   E   M   A   R   K      %%%%%%%%%%%%%%%
%%%%%%%%%%%%%%%%%%%%%%%%%%%%%%%%
%%%%%%%%%%%%%%%%%%%%%%%%%%%%%%%% 
\begin{remark}
One may raise the question about the choice of $K$ for dynamic ESOM-$K$. Note that the implementation of dynamic ESOM-$K$ requires $K+1$ rounds of communication between neighboring nodes. Thus, by increasing the choice of $K$ the computation time of the algorithm increases. Although, larger choice of $K$ leads to a better Hessian inverse approximation and faster convergence, the required time may exceed the time between the subsequent instances $t-1$ and $t$. Therefore, based on the available time between the consecutive times $t-1$ and $t$, we should pick the largest choice of $K$ which is affordable in terms of computation and communication time.
\end{remark}

%%%%%%%%%%%%%%%%%%%%%%%%%%%%%%%%
%%%%%%%%%%%%%%%%%%%%%%%%%%%%%%%%
%%%%%     S   E   C   T   I   O   N      %%%%%%%%%%%%%
%%%%%%%%%%%%%%%%%%%%%%%%%%%%%%%%
%%%%%%%%%%%%%%%%%%%%%%%%%%%%%%%%
\section{Convergence Analysis}\label{sec_convg}
In this section we study the difference between the sequence of the iterates $\bbx_t$ generated by dynamic ESOM and the sequence of the optimal arguments $\bbx_t^*=[\bbx_{1,t}^*;\dots;\bbx_{n,t}^*]=[\tbx_t^*;\dots;\tbx_t^*]$. To prove the results, we assume the following conditions are satisfied.

%%%%%%%%%%%%%%%%%%%%%%%%%%%%%%%%
%%%%%%%%%%%%%%%%%%%%%%%%%%%%%%%%
%%%    A   S   S   U   M   P   T   I   O   N    %%%%%%%%%%%
%%%%%%%%%%%%%%%%%%%%%%%%%%%%%%%%
%%%%%%%%%%%%%%%%%%%%%%%%%%%%%%%%
\begin{assumption}\label{convexity_assumption} 
The instantaneous local objective functions $f_{i,t}(\bbx)$ are twice differentiable and the eigenvalues of the instantaneous local objective functions Hessian $\nabla^2 f_{i,t}$ are bounded by positive constants $0<m\leq M<\infty$, i.e. 
\begin{equation}\label{local_hessian_eigenvlaue_bounds00}
m\bbI  \ \preceq\ \nabla^2 f_{i,t}(\bbx_i)\ \preceq\ M\bbI, 
\end{equation}
for all $\bbx_i \in \reals^p$ and $i=1,\dots,n$.
\end{assumption}

%%%%%%%%%%%%%%%%%%%%%%%%%%%%%%%%
%%%%%%%%%%%%%%%%%%%%%%%%%%%%%%%%
%%%    A   S   S   U   M   P   T   I   O   N    %%%%%%%%%%%
%%%%%%%%%%%%%%%%%%%%%%%%%%%%%%%%
%%%%%%%%%%%%%%%%%%%%%%%%%%%%%%%%
%
\begin{assumption}\label{lip_assumption} 
The instantaneous local objective functions Hessian $\nabla^2f_{i,t}$ are Lipschitz continuous with constant $L$, 
\begin{equation}\label{local_hessian_eigenvlaue_bounds100}
\| \nabla^2 f_{i,t}(\bbx_i)- \nabla^2 f_{i,t}(\bby_i)\| \leq L\|\bbx_i-\bby_i\|, 
\end{equation}
for all $\bbx_i,\bby_i \in \reals^p$ and $i=1,\dots,n$.
\end{assumption}

%%%%%%%%%%%%%%%%%%%%%%%%%%%%%%%%
%%%%%%%%%%%%%%%%%%%%%%%%%%%%%%%%
%%%         M  A  I  N        M  A  T  T  E  R     %%%%%%%%%%
%%%%%%%%%%%%%%%%%%%%%%%%%%%%%%%%
%%%%%%%%%%%%%%%%%%%%%%%%%%%%%%%% 

We can interpret the lower and upper bounds on the eigenvalues of the Hessians $\nabla^2f_{i,t}$ as the strong convexity of the instantaneous local functions $f_{i,t}$ with constant $m$ and the Lipschitz continuity of the instantaneous local gradients $\nabla f_{i,t}$ with constant $M$, respectively. The global objective function Hessian $\nabla^2 f_t(\bbx)$ at step $t$ is a block diagonal matrix where its $i$-th diagonal block is $\nabla^2 f_{i,t}(\bbx_i)$. Hence, the bounds in \eqref{local_hessian_eigenvlaue_bounds00} for the eigenvalues of the instantaneous local Hessians also hold for the instantaneous global Hessian $\nabla^2 f_t(\bbx)$, i.e.,
\begin{equation}\label{local_hessian_eigenvlaue_bounds}
m\bbI \ \preceq \ \nabla^2 f_t(\bbx) \ \preceq\  M\bbI,
\end{equation}
for all $\bbx \in \reals^{np}$. Thus, the global objective function $f_t$ is also strongly convex with constant $m$ and its gradients $\nabla f_t$ are Lipschitz continuous with constant $M$. Likewise, the Lipschitz continuity of the local Hessians $\nabla^2f_{i,t}$, which is a customary assumption in the analysis of second-order methods, implies that the instantaneous global Hessian $ \nabla^2 f_{t}$ is also Lipschitz continuous with constant $L$, i.e., 
\begin{equation}\label{local_hessian_eigenvlaue_bounds200}
\| \nabla^2 f_{t}(\bbx)- \nabla^2 f_{t}(\bby)\| \leq L\|\bbx-\bby\|, 
\end{equation}
for any $\bbx,\bby\in \reals^{np}$ -- see e.g., Lemma 1 in \cite{NN-part1}. 

To characterize the error of dynamic ESOM, we define the vector $\bbu_t=[\bbx_t; \bbv_t]\in \reals^{2np}$ as the concatenation of the primal and dual iterates at step $t$. Likewise, we define $\bbu_t^*=[\bbx_t^*; \bbv_t^*]\in \reals^{2np}$ as the concatenation of the optimal arguments at time $t$. We proceed to characterize an upper bound for the error sequence $\|\bbu_t-\bbu_t^*\|_\bbG$ where the positive definite matrix $\bbG$ is defined as $\bbG:=\diag(\bbI_{np},\eps \alpha \bbI_{np})\in \reals^{2np\times2np}$. In the following lemma we establish an upper bound for the norm $\|\bbu_t-\bbu_t^*\|_\bbG$ in terms of the difference between the previous vector $\bbu_{t-1}$ and the current optimal argument $\bbu_t^*$.

%%%%%%%%%%%%%%%%%%%%%%%%%%%%%%%%
%%%%%%%%%%%%%%%%%%%%%%%%%%%%%%%%
%%%%%        L   E   M   M    A        %%%%%%%%%%%%%%
%%%%%%%%%%%%%%%%%%%%%%%%%%%%%%%%
%%%%%%%%%%%%%%%%%%%%%%%%%%%%%%%%
\begin{lemma}\label{lemma:pmm_linear_convg}
Consider the updates of dynamic ESOM as introduced in \eqref{ESOM_primal_update_3}-\eqref{ESOM_dual_update} and recall the definitions of the vector $\bbu$ and matrix $\bbG$. If Assumptions \ref{convexity_assumption} and \ref{lip_assumption} hold, then there exists a positive scalar $0<\delta$ such that the sequence of iterates $\bbu_t$ generated by dynamic ESOM satisfies 
\begin{equation}\label{proof_0}
\|\bbu_{t}-\bbu_t^*\|_\bbG \ \leq\  \frac{1}{\sqrt{1+\delta}} \ \|\bbu_{t-1}-\bbu_t^*\|_\bbG.
\end{equation}
%%%
\end{lemma}

%%%%%%%%%%%%%%%%%%%%%%%%%%%%%%%%
%%%%%%%%%%%%%%%%%%%%%%%%%%%%%%%%
%%%%%         P    R    O    O     F      %%%%%%%%%%%%%
%%%%%%%%%%%%%%%%%%%%%%%%%%%%%%%%
%%%%%%%%%%%%%%%%%%%%%%%%%%%%%%%%
\begin{myproof}
The proof can be established by following the steps of the proof of Theorem 2 in \cite{mokhtari2016decentralized}.
\end{myproof}

%%%%%%%%%%%%%%%%%%%%%%%%%%%%%%%%
%%%%%%%%%%%%%%%%%%%%%%%%%%%%%%%%
%%%         M  A  I  N        M  A  T  T  E  R     %%%%%%%%%%
%%%%%%%%%%%%%%%%%%%%%%%%%%%%%%%%
%%%%%%%%%%%%%%%%%%%%%%%%%%%%%%%% 

The constant $\delta$ in \eqref{proof_0} is a function of the objective function $f_t$ parameters, network topology, and level of Hessian inverse approximation $K$. In particular, the constant $\delta$ is close to zero when the objective function is ill-conditioned, or the network is not well connected. Moreover, larger choice of $K$ leads to a larger choice of $\delta$ which leads to a smaller error $\|\bbu_{t}-\bbu_t^*\|_\bbG$.

The result in Lemma \ref{lemma:pmm_linear_convg} illustrates that the iterate $\bbu_t$ is closer to the optimal argument $\bbu_t^*$ at step $t$ relative to the previous iterate $\bbu_{t-1}$. This result is implied from the fact that $\bbu_t$ is evaluated based on the observed function $f_t$ at step $t$. 
Based on the result in Lemma \ref{lemma:pmm_linear_convg}, we can establish an upper bound for the error $\|\bbu_{t}-\bbu_t^*\|_\bbG$ at step $t$ in terms of the error of the previous time $\|\bbu_{t-1}-\bbu_{t-1}^*\|_\bbG$ and the variation of the optimal arguments. We characterize this upper bound in the following theorem.

%%%%%%%%%%%%%%%%%%%%%%%%%%%%%%%%
%%%%%%%%%%%%%%%%%%%%%%%%%%%%%%%%
%%%%%        T   H   E   O   R   E   M        %%%%%%%%%%%
%%%%%%%%%%%%%%%%%%%%%%%%%%%%%%%%
%%%%%%%%%%%%%%%%%%%%%%%%%%%%%%%%
\begin{thm}\label{lin_thm}
Consider the dynamic ESOM algorithm as introduced in \eqref{ESOM_primal_update_3}-\eqref{ESOM_dual_update} and recall the definitions of the vector $\bbu$ and matrix $\bbG$. Define $\gamma$ as the smallest non-zero eigenvalue of the positive semidefinite matrix $\bbI-\bbZ$. Further, define the dynamic optimality drift $d_t$ as
\begin{equation}\label{drift}
d_t :=  \|\bbx_{t-1}^*-\bbx_t^*\|+\frac{\sqrt{\alpha \eps} }{\sqrt{\gamma}} \|\nabla f_t(\bbx_t^*)-\nabla f_{t-1}(\bbx_{t-1}^*) \|.
\end{equation}
If Assumptions \ref{convexity_assumption} and \ref{lip_assumption} hold, then the sequence of iterates $\bbu_t$ generated by dynamic ESOM satisfies 
%%%
\begin{equation}\label{proof_011}
\|\bbu_{t}-\bbu_t^*\|_\bbG \ \leq\ \frac{1}{\sqrt{1+\delta}} \ \! \|\bbu_{t-1}-\bbu_{t-1}^*\|_\bbG +\frac{d_t}{\sqrt{1+\delta}}.
\end{equation}
\end{thm}

%%%%%%%%%%%%%%%%%%%%%%%%%%%%%%%%
%%%%%%%%%%%%%%%%%%%%%%%%%%%%%%%%
%%%%%         P    R    O    O     F      %%%%%%%%%%%%%
%%%%%%%%%%%%%%%%%%%%%%%%%%%%%%%%
%%%%%%%%%%%%%%%%%%%%%%%%%%%%%%%%
\begin{myproof}
See Appendix \ref{lemma:new_lemma_app}.
\end{myproof}

%%%%%%%%%%%%%%%%%%%%%%%%%%%%%%%%
%%%%%%%%%%%%%%%%%%%%%%%%%%%%%%%%
%%%         M  A  I  N        M  A  T  T  E  R     %%%%%%%%%%
%%%%%%%%%%%%%%%%%%%%%%%%%%%%%%%%
%%%%%%%%%%%%%%%%%%%%%%%%%%%%%%%% 
The optimality drift $d_t$ captures the drift between the two consecutive optimal arguments $\bbx_t^*$ and $\bbx_{t-1}^*$ as well as the difference between the two successive optimal gradients $\nabla f_t(\bbx_t^*)$ and $\nabla f_{t-1}(\bbx_{t-1}^*)$. The result in Theorem \ref{lin_thm} shows that the sequence of the error $\|\bbu_{t}-\bbu_t^*\|_\bbG$ approaches linearly a steady state error bound. Note that the optimality drift $d_t$ is small when the functions $f_t$ change sufficiently slow. The result in \eqref{proof_011} is consistent with the results for the static version of the optimization problem in \eqref{gen_prob}. In the static setting, where $\bbu_t^*=\bbu_{t-1}^*=\bbu^*$, $\bbx_t^*=\bbx_{t-1}^*=\bbx^*$, and $\nabla f_t(\bbx_t^*)=\nabla f_{t-1}(\bbx_{t-1}^*)=\nabla f(\bbx^*)$, the result in \eqref{proof_011} can be simplified as $\|\bbu_{t}-\bbu^*\|_\bbG  \leq ({1}/{\sqrt{1+\delta}}) \|\bbu_{t-1}-\bbu^*\|_\bbG$ which shows linear convergence of the iterates generated by ESOM to the optimal argument.

In the following theorem we use the result in Theorem \ref{lin_thm} to show that the error of dynamic ESOM, which is characterized by the norm $\|\bbu_{t}-\bbu_t^*\|_\bbG$, approaches a steady state error.

%%%%%%%%%%%%%%%%%%%%%%%%%%%%%%%%
%%%%%%%%%%%%%%%%%%%%%%%%%%%%%%%%
%%%%%        T   H   E   O   R   E   M        %%%%%%%%%%%
%%%%%%%%%%%%%%%%%%%%%%%%%%%%%%%%
%%%%%%%%%%%%%%%%%%%%%%%%%%%%%%%%
\begin{thm}\label{theorem:new_theorem}
Consider the dynamic ESOM algorithm as introduced in \eqref{ESOM_primal_update_3}-\eqref{ESOM_dual_update} and recall the definition of the optimality drift $d_t$ in \eqref{drift}. Further, define $d_{\max}:=\max_{t}d_t$ as the maximum of the optimality drift $d_t$ for all times $t$.
If Assumptions \ref{convexity_assumption} and \ref{lip_assumption} hold, then the limit supremum of the sequence $\|\bbu_{t}-\bbu_t^*\|_\bbG$ is bounded above by
%%%
\begin{align}\label{claim_thm_upper_bound_sup}
\limsup_{t\to \infty} \|\bbu_{t}-\bbu_t^*\|_\bbG 
\ \leq\   \frac{d_{\max}}{\sqrt{1+\delta}-1}.
\end{align}
\end{thm}

%%%%%%%%%%%%%%%%%%%%%%%%%%%%%%%%
%%%%%%%%%%%%%%%%%%%%%%%%%%%%%%%%
%%%%%         P    R    O    O     F      %%%%%%%%%%%%%
%%%%%%%%%%%%%%%%%%%%%%%%%%%%%%%%
%%%%%%%%%%%%%%%%%%%%%%%%%%%%%%%%
\begin{myproof}
See Appendix \ref{theorem:new_theorem_app}.
\end{myproof}

%%%%%%%%%%%%%%%%%%%%%%%%%%%%%%%%
%%%%%%%%%%%%%%%%%%%%%%%%%%%%%%%%
%%%         M  A  I  N        M  A  T  T  E  R     %%%%%%%%%%
%%%%%%%%%%%%%%%%%%%%%%%%%%%%%%%%
%%%%%%%%%%%%%%%%%%%%%%%%%%%%%%%% 
The steady state error of the sequence generated by dynamic ESOM is characterized in Theorem \ref{theorem:new_theorem}. As we expect, if the maximum optimality drift $d_{\max}$ is not  large the dynamic ESOM algorithm approaches a reasonable asymptotic error. Moreover, the steady state error is smaller for the case that the constant of linear convergence $\delta$ is larger. This observation shows that the steady state error of ESOM-$K$ reduces by increasing the level of Hessian inverse approximation $K$. This is true, since for larger choice of $K$, the constant $\delta$ is larger. 

Convergence of the sequence $\|\bbu_{t}-\bbu_t^*\|_\bbG$, which characterizes the primal and dual errors of the iterates of dynamic ESOM, follows that the sequence of primal iterates $\bbx_t$ converges to a neighborhood of the optimal argument $\bbx_t^*$. This is shown in the following corollary.

%%%%%%%%%%%%%%%%%%%%%%%%%%%%%%%%
%%%%%%%%%%%%%%%%%%%%%%%%%%%%%%%%
%%%%%         C  O  R  O  L  L   A   R  Y      %%%%%%%%%%
%%%%%%%%%%%%%%%%%%%%%%%%%%%%%%%%
%%%%%%%%%%%%%%%%%%%%%%%%%%%%%%%%
\begin{corollary}\label{corollary1}
Recall the definition of the maximum optimality drift $d_{\max}$ and suppose that the conditions in Theorem \ref{theorem:new_theorem} are satisfied. Then the primal error $\|\bbx_t-\bbx_t^*\|$ of dynamic ESOM is upper bounded as
%%%
\begin{align}\label{claim_cor_upper_bound_sup}
\limsup_{t\to \infty} \|\bbx_{t}-\bbx_t^*\| 
\ \leq\   \frac{d_{\max}}{\sqrt{1+\delta}-1}.
\end{align}
%%%
\end{corollary}

%%%%%%%%%%%%%%%%%%%%%%%%%%%%%%%%
%%%%%%%%%%%%%%%%%%%%%%%%%%%%%%%%
%%%%%         P    R    O    O     F      %%%%%%%%%%%%%
%%%%%%%%%%%%%%%%%%%%%%%%%%%%%%%%
%%%%%%%%%%%%%%%%%%%%%%%%%%%%%%%%
\begin{myproof}
Based on the definition of the norm $\|\bbu_{t}-\bbu_t^*\|_\bbG$, we can simplify the norm as $\|\bbu_{t}-\bbu_t^*\|_\bbG=\left[\|\bbx_{t}-\bbx_t^*\|^2+\alpha \eps\|\bbv_{t}-\bbv_t^*\|^2\right]^{1/2}$. According to this defnition, we obtain that the primal error $\|\bbx_{t}-\bbx_t^*\|$ is smaller than the norm $\|\bbu_{t}-\bbu_t^*\|_\bbG$. This observation in conjunction with the result in \eqref{claim_thm_upper_bound_sup} implies the claim in \eqref{claim_cor_upper_bound_sup}.
\end{myproof}

%%%%%%%%%%%%%%%%%%%%%%%%%%%%%%%%
%%%%%%%%%%%%%%%%%%%%%%%%%%%%%%%%
%%%%%     S   E   C   T   I   O   N      %%%%%%%%%%%%%
%%%%%%%%%%%%%%%%%%%%%%%%%%%%%%%%
%%%%%%%%%%%%%%%%%%%%%%%%%%%%%%%%
\section{Numerical experiments}\label{sec:num_exp}

In this section, we study the performance of the proposed dynamic ESOM method in solving a dynamic least squares problem. We consider a connected network with $n=20$ and connectivity ratio $r_c=0.15$, i.e., edges are generated randomly with probability $0.15$.

We consider a decentralized dynamic least squares problem where at time $t$ nodes aim to estimate the true signal $\tbx^*_t\in \reals^{5}$. Consider the linear model $\bby_{i,t}=\bbH_{i,t}\tbx^*_t+\bbeta_{i,t}$ where the matrix $\bbH_{i,t}\in \reals^{5\times 5}$ is a regressor matrix and the vector $\bbeta_{i,t}\in \reals^{5}$ is an additive noise. We assume that node $i$ observes the vector $\bby_{i,t}$ and collaborates with its neighbors to find the true signal $\tbx_t^*$. In other words, the nodes' goal is to solve the least squares problem 
\begin{equation}\label{cost_simulation}
\tbx_t^*=\argmin_{\bbx\in \reals^5} \sum_{i=1}^n \frac{1}{2}\|\bbH_{i,t}\bbx-\bby_{i,t}\|^2.
\end{equation}
Considering the definition of the global optimization problem in \eqref{cost_simulation}, the local objective function of node $i$ at time $t$ is given by $f_{i,t}:=(1/2)\|\bbH_{i,t}\bbx-\bby_{i,t}\|^2$.

We compare the dynamic variations of ESOM-$0$, ESOM-$2$, Network Newton-0 (NN-$0$) \cite{NN-part1}, and EXTRA \cite{Shi2014} in solving the dynamic least squares problem in \eqref{cost_simulation}. In our experiments, we assume that the matrices $\bbH_{i,t} $ are fixed over time, i.e., $\bbH_{i,t}=\bbH_{i}$. We generate the components of $\bbH_{i}$ following the Gaussian distribution $\ccalN(0,1)$. Although, the matrices $\bbH_{i}$ are time-invariant, we assume that the vectors $\bby_{i,t}$ are changing over time. We assume that after every $100$ iterations the components of the vectors $\bby_{i,t}$ change in a way that the new global minimizer $\tbx_t^*$ satisfies $\tbx_t^*=|\sin(\pi t/500)|\tbx_{0}^*$. In other words, $\tbx_t^*=|\sin(\pi t/500)|\tbx_{0}^*$ if $t$ is a multiplicant of 100, otherwise  $\tbx_t^*=\tbx_{t-1}^*$. Moreover, we assume that every agent starts from the initial point $\bbx_{i,0}$ that satisfies the condition $\|\bbx_{i,0}-\tbx_{0}^*\|=100$.

We characterize error as the maximum difference between the coordinates of each node's variable and the optimal argument $\tbx_t^*$. Thus, if we define $\bbx_{i,t}[s]$ as the $s$-th coordinate of the variable $\bbx_{i,t}$, the error is defined as $e_t:=\max_{i}\{\max_{s}\{|\bbx_{i,t}[s]-\tbx_{t}^*[s]|\}\}$. The error $e_t$ also can be  written as
\begin{equation}
e_t:=\max_{i}\left\{\|\bbx_{i,t}-\tbx_{t}^*\|_{\infty} \right\},
\end{equation}
using the definition of the infinity norm $\|\cdot\|_{\infty}$.

%%%%%%%%%%% FIGURE %%%%%%%%%%%%%%%%%
\begin{figure}[t]
\centering
\includegraphics[width=\linewidth]{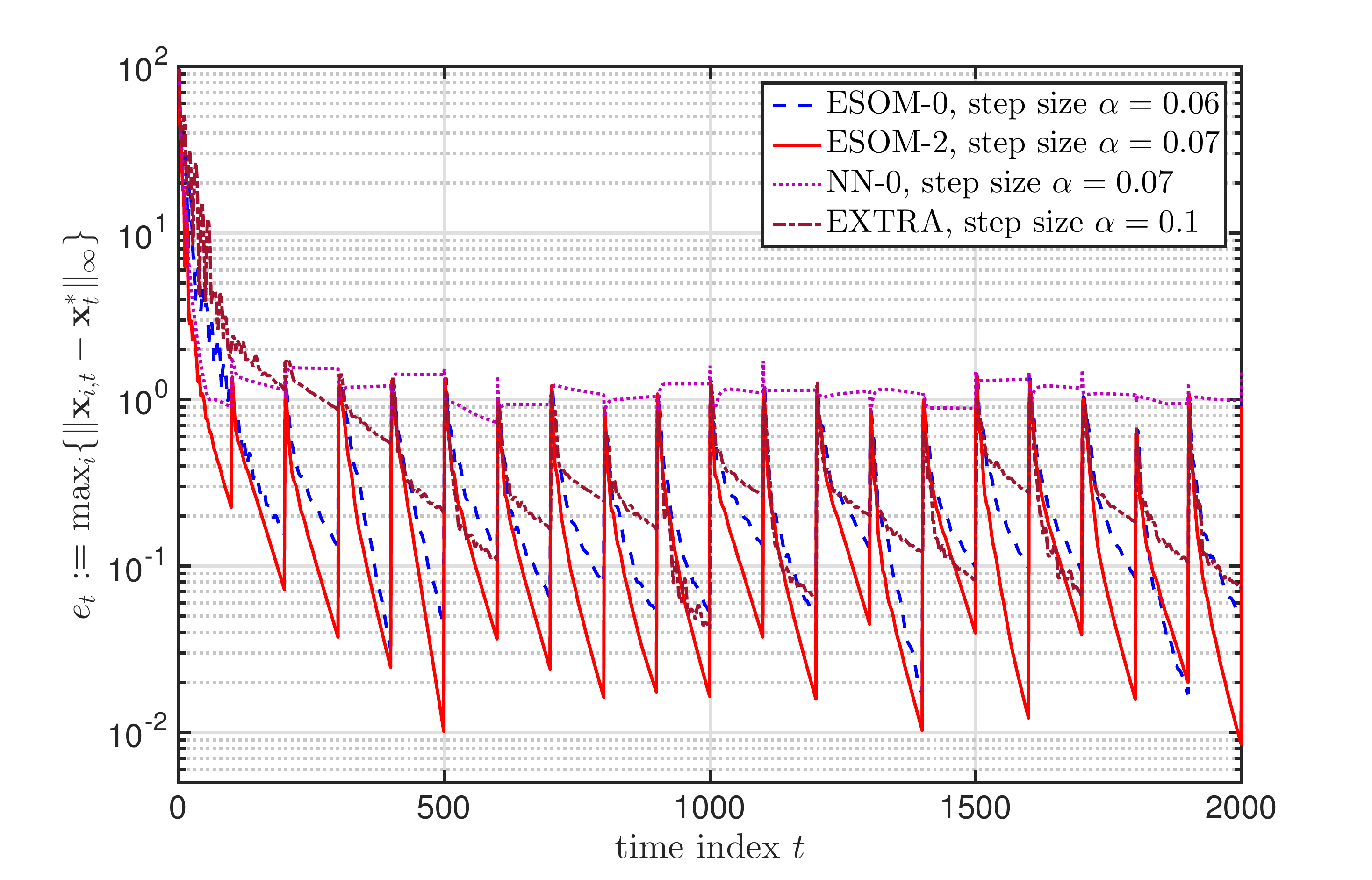}
\caption{The convergence path of $e_t$ versus time index $t$ for dynamic ESOM-$0$, ESOM-$2$, EXTRA, and NN-$0$. Dynamic ESOM-$2$ has the best performance among all the considered methods.}\label{figure_simulation1}
\end{figure}
%%%%%%%%%%%%%%%%%%%%%%%%%%%%%%%%%

Figure \ref{figure_simulation1} shows the error $e_t$ versus the time index $t$ for the four algorithms of interest. As we observe, during the time that the optimal argument is fixed, NN-$0$ approaches a neighborhood of the optimal solution and its error $e_t$ stays constant, while EXTRA, ESOM-$0$, and ESOM-$2$ converge linearly to the exact solution and their error $e_t$ diminish. It is also worth mentioning that both ESOM-$0$ and ESOM-$2$ outperform EXTRA by incorporating second-order information, and ESOM-$2$ has the best performance among all the considered methods. If more rounds of communication is affordable between the subsequent instances $t$ and $t+1$, then the performance of dynamic ESOM-$K$ can be improved by using larger values for $K$. Note that whenever the optimal argument $\bbx_t^*$ changes, which happens every $100$ iterations, all the algorithms readjust and correct their descent direction to track the new optimal argument. 

\begin{figure}[t]
\centering
\includegraphics[width=\linewidth]{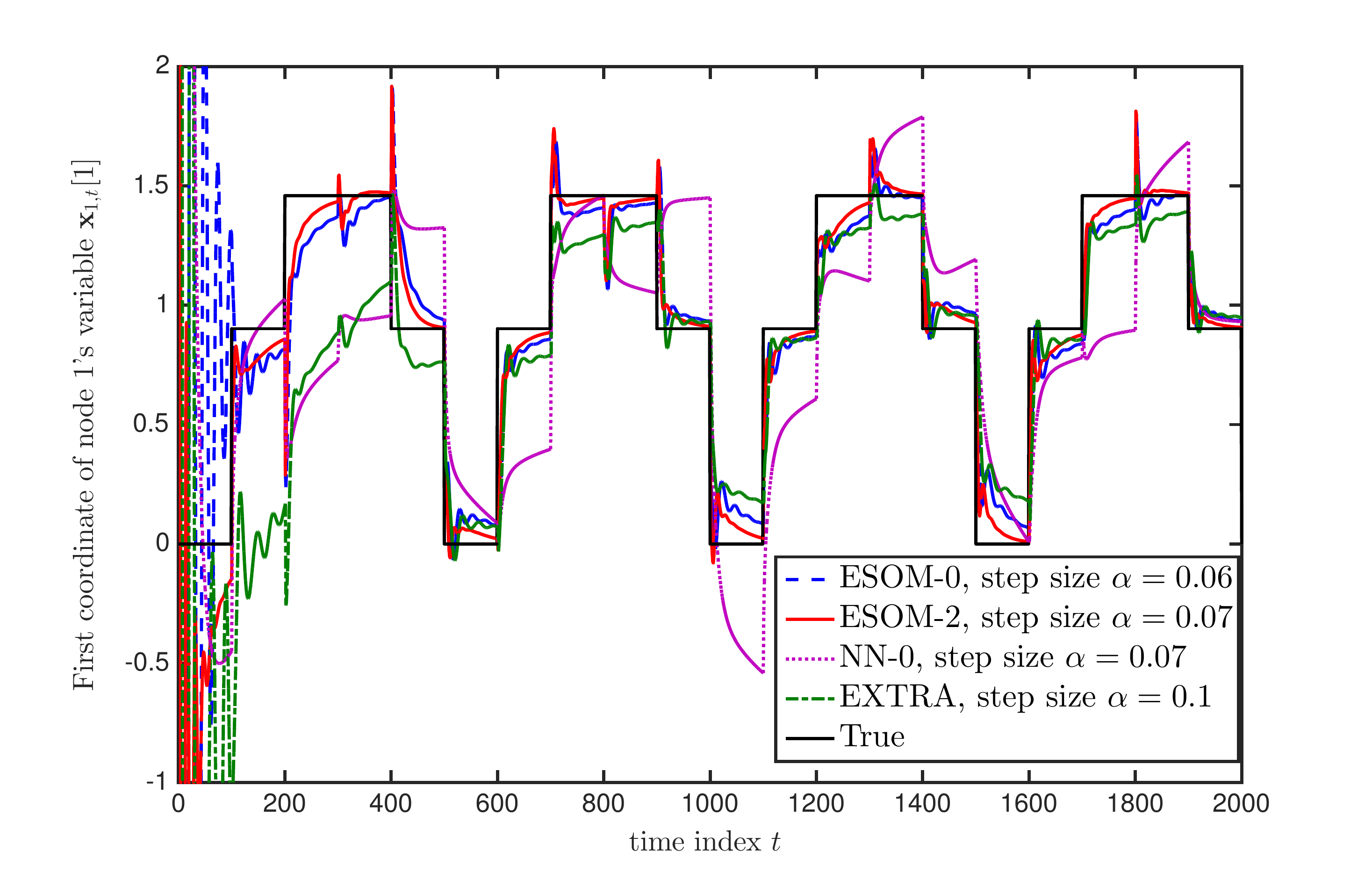}
\caption{Comparison of node 1's first coordinate for the iterates generated by dynamic variations of ESOM-$0$, ESOM-$2$, NN-$0$, and EXTRA with the first coordinate of the optimal argument $\tbx_{t}^*$.}\label{figure_simulation2}
\end{figure}

To study the performance of dynamic NN-$0$, EXTRA, ESOM-$0$, and ESOM-$2$ in more details, we compare the values of the first coordinate $\bbx_{1,t}[1]$ of node 1 generated by these methods with the first coordinate of the optimal argument $\tbx_t^*[1]$. This comparison is shown in Figure \ref{figure_simulation2}. As we observe in Figure \ref{figure_simulation2}, all the dynamic methods are unable to track the true path in the first $200$ iterations. Dynamic ESOM-$0$ and ESOM-$2$ can track the optimal argument after the first $200$ iterations, while the accuracy of dynamic ESOM-$2$ is higher relative to dynamic ESOM-$0$. Dynamic EXTRA starts tracking the true path after $400$ iterations, while its error is worse than the ones for dynamic ESOM-$0$ and ESOM-$2$. For the dynamic NN-$0$ method, the error is always larger than the error of the other dynamic methods. 

These observations verify the theoretical results in Section \ref{sec_convg}. To be more precise, they show that the dynamic variations of EXTRA and ESOM-$K$ outperform dynamic NN, since they converge linearly to the optimal arguments while NN converges to a neighborhood of the optimal solution in static settings. Moreover, dynamic ESOM-$K$, irrespective to the choice of $K$, improves the performance of dynamic EXTRA by incorporating second-order information of the augmented Lagrangian in \eqref{lagrangian}. Further, larger choice of $K$ for dynamic ESOM-$K$ leads to a faster linear convergence, i.e., larger $\delta$, which implies a smaller steady state error. This observation verifies the result in Corollary \ref{corollary1}. 

%%%%%%%%%%% FIGURE %%%%%%%%%%%%%%%%%
%%%%%%%%%%%%%%%%%%%%%%%%%%%%%%%%%

%%%%%%%%%%%%%%%%%%%%%%%%%%%%%%%%
%%%%%%%%%%%%%%%%%%%%%%%%%%%%%%%%
%%%%%     S   E   C   T   I   O   N      %%%%%%%%%%%%%
%%%%%%%%%%%%%%%%%%%%%%%%%%%%%%%%
%%%%%%%%%%%%%%%%%%%%%%%%%%%%%%%%
\section{Conclusions}\label{sec:conclusions}
{We considered the application of the Exact Second-Order Methods (ESOM) in solving a dynamic consensus optimization problem where the local functions available at nodes are time-variant. The proposed dynamic ESOM method relies on the use of a separable quadratic approximation of a suitably defined time-varying augmented Lagrangian, and a truncated Taylor's series to estimate the solution of the first order condition imposed on the minimization of the quadratic approximation
of the augmented Lagrangian. We proved that under proper assumptions, the sequence of iterates generated by dynamic ESOM converges linearly to a neighborhood of the sequence of optimal arguments. We characterized the steady state error in terms of the maximum difference between the successive optimal arguments $\bbx_{t-1}^*$ and $\bbx_{t}^*$ as well as the optimal gradients $\nabla f_{t-1}(\bbx_{t-1}^*)$ and $\nabla f_{t}(\bbx_{t}^*) $. Numerical results showcase the advantages of the proposed dynamic ESOM method relative to existing dynamic decentralized methods.

%%%%%%%%%%%%%%%%%%%%%%%%%%%%%%%%
%%%%%%%%%%%%%%%%%%%%%%%%%%%%%%%%
%%%%%     S   E   C   T   I   O   N      %%%%%%%%%%%%%
%%%%%%%%%%%%%%%%%%%%%%%%%%%%%%%%
%%%%%%%%%%%%%%%%%%%%%%%%%%%%%%%%
\section{APPENDIX}

\subsection{Proof of Theorem \ref{lin_thm}}\label{lemma:new_lemma_app}

According to the definition of the vector $\bbu$ and matrix $\bbG$ we can write 
%%%
\begin{align}\label{jadidan}
\|\bbu_{t-1}^*-\bbu_t^*\|_\bbG
&=\left[\|\bbx_{t-1}^*-\bbx_t^*\|^2+\alpha \eps\|\bbv_{t-1}^*-\bbv_t^*\|^2\right]^{1/2}\nonumber\\
&\leq \|\bbx_{t-1}^*-\bbx_t^*\|+\sqrt{\alpha \eps}\|\bbv_{t-1}^*-\bbv_t^*\|,
\end{align}
where the inequality follows from the inequality $a^2+b^2\leq (a+b)^2$ for positive scalars $a$ and $b$. The KKT condition of the optimization problem in \eqref{constrained_opt_problem} yields
\begin{equation}\label{KKT_condition_22}
\nabla f_t(\bbx_t^*) + (\bbI-\bbZ)^{1/2}\bbv_t^*=\bb0.
\end{equation}
%%%
By writing the KKT condition in \eqref{KKT_condition_22} for time $t-1$ we obtain that 
\begin{equation}\label{KKT_condition_223}
\nabla f_t(\bbx_t^*)-\nabla f_{t-1}(\bbx_{t-1}^*) + (\bbI-\bbZ)^{1/2}(\bbv_t^*-\bbv_{t-1}^*)=\bb0.
\end{equation}
%%%
Since the vectors $\bbv_t^*$ and $\bbv_{t-1}^*$ are in the column space of the matrix $\bbI-\bbZ$, we obtain that $ (\bbI-\bbZ)^{1/2}(\bbv_t^*-\bbv_{t-1}^*)$ is bounded below by $\gamma^{1/2}\|\bbv_t^*-\bbv_{t-1}^*\|$. This lower bound in conjunction with the expression in \eqref{KKT_condition_223} implies that the norm $\|\bbv_t^*-\bbv_{t-1}^*\|$ is bounded above as
\begin{equation}\label{KKT_condition_2123}
\|\bbv_t^*-\bbv_{t-1}^*\| \leq \frac{1}{\sqrt{\gamma}} \|\nabla f_t(\bbx_t^*)-\nabla f_{t-1}(\bbx_{t-1}^*) \|.
\end{equation}
%%%
%%
By substituting the upper bound in \eqref{KKT_condition_2123} to \eqref{jadidan}, the inequality $\|\bbu_{t-1}^*-\bbu_t^*\|_\bbG\leq d_t$ follows.

Based on the triangle inequality, the weighted norm $\|\bbu_{t-1}-\bbu_t^*\|_\bbG$ is bounded above by the sum
\begin{equation}\label{khosro}
\|\bbu_{t-1}-\bbu_t^*\|_\bbG\leq  \|\bbu_{t-1}-\bbu_{t-1}^*\|_\bbG+\|\bbu_{t-1}^*-\bbu_t^*\|_\bbG.
 \end{equation}
Since the norm $\|\bbu_{t-1}^*-\bbu_t^*\|_\bbG$ is smaller than the drift $d_t$ defined in \eqref{drift}, we can replace $\|\bbu_{t-1}^*-\bbu_t^*\|_\bbG$ in \eqref{khosro} by $d_t$
%%%
\begin{equation}\label{proof_0989811}
\|\bbu_{t-1}-\bbu_t^*\|_\bbG \ \! \leq \ \! \|\bbu_{t-1}-\bbu_{t-1}^*\|_\bbG +d_t.
\end{equation}
Combining the results in \eqref{proof_0989811} and \eqref{proof_0}, the claim in \eqref{proof_011} follows.

%%%%%%%%%%%%%%%%%%%%%%%%%%%%%%%%%%%%%%%%%%%%%%%%%%%%%%%%%%%%%%%%%%
%%%   APPENDIX  -- C  %%%%%%%%%%%%%%%%%%%%%%%%%%%%%%%%%%%%%%%%
%%%%%%%%%%%%%%%%%%%%%%%%%%%%%%%%%%%%%%%%%%%%%%%%%%%%%%%%%%%%%%%%%%
%
\subsection{Proof of Theorem \ref{theorem:new_theorem}}\label{theorem:new_theorem_app}

We prove the claim in \eqref{claim_thm_upper_bound_sup} based on \eqref{proof_011} in Theorem \ref{lin_thm}. First, consider the inequality in \eqref{proof_011} for time $t-1$ which is given by 
\begin{equation}\label{proof_0921}
\|\bbu_{t-1}-\bbu_{t-1}^*\|_\bbG \ \leq\ \frac{1}{\sqrt{1+\delta}} \ \! \|\bbu_{t-2}-\bbu_{t-2}^*\|_\bbG +\frac{d_{t-1}}{\sqrt{1+\delta}}.
\end{equation}
Substituting the norm $\|\bbu_{t-1}-\bbu_{t-1}^*\|_\bbG$ in \eqref{proof_011} by the upper bound in \eqref{proof_0921} yields
%%%
\begin{equation}\label{proof_new_proof_10949}
\|\bbu_{t}-\bbu_t^*\|_\bbG \ \leq\ \! \frac{\|\bbu_{t-2}-\bbu_{t-2}^*\|_\bbG}{(\sqrt{1+\delta})^2} + \frac{d_{t-1}}{(\sqrt{1+\delta})^2}+\frac{d_t}{\sqrt{1+\delta}}.
\end{equation}
By considering the expression in \eqref{proof_011} for all times $s\leq t$ and recursively it follows that the error $ \|\bbu_{t}-\bbu_t^*\|_\bbG$ at step $t$ and the initial error $ \|\bbu_{0}-\bbu_0^*\|_\bbG$ satisfy
%%%
\begin{equation}\label{proof_new_proof_100}
\|\bbu_{t}-\bbu_t^*\|_\bbG \ \leq\ \! \frac{\|\bbu_{0}-\bbu_{0}^*\|_\bbG}{(\sqrt{1+\delta})^t} +\sum_{s=1}^t \frac{d_s}{(\sqrt{1+\delta})^{t-s+1}}.
\end{equation}
Now considering the definition of $d_{\max}$ as $d_{\max}=\sup_{t\geq1}d_t$, we obtain that the sum in the right hand side of \eqref{proof_new_proof_100} is bounded above by
%%%
\begin{align}\label{proof_new_proof_1900}
\sum_{s=1}^t \frac{d_s}{(\sqrt{1+\delta})^{t-s+1}}& \ \leq\ \sum_{s=1}^t \frac{d_{\max}}{(\sqrt{1+\delta})^{t-s+1}}\nonumber\\
& \ \leq\ \frac{d_{\max}}{\sqrt{1+\delta}} \times  \frac{1-(\sqrt{1+\delta})^{-t}}{1-(\sqrt{1+\delta})^{-1}},
\end{align}
where the second inequality is implied from the simplification $\sum_{s=1}^t \rho^t=\rho (1-\rho^t)/(1-\rho)$ when $\rho<1$. Replacing the sum in the right hand side of \eqref{proof_new_proof_100} by the upper bound in \eqref{proof_new_proof_1900} yields 
%%%
\begin{equation}\label{proof_new_proof_1090}
\|\bbu_{t}-\bbu_t^*\|_\bbG \ \leq\ \! \frac{\|\bbu_{0}-\bbu_{0}^*\|_\bbG}{(\sqrt{1+\delta})^t} +\frac{d_{\max}}{\sqrt{1+\delta}} \times  \frac{1-(\sqrt{1+\delta})^{-t}}{1-(\sqrt{1+\delta})^{-1}}.
\end{equation}
Taking $t\to\infty$, the term ${\|\bbu_{0}-\bbu_{0}^*\|_\bbG}/{(\sqrt{1+\delta})^t}$ in the right hand side of \eqref{proof_new_proof_1090} vanishes. Moreover, the term $(\sqrt{1+\delta})^{-t}$ approaches $0$. From these observations it follows that 
%%%
\begin{align}\label{proof_new_proof_200}
\limsup_{t\to \infty} \|\bbu_{t}-\bbu_t^*\|_\bbG 
\ \leq\   \frac{d_{\max}}{\sqrt{1+\delta}-1},
\end{align}
and the proof is complete.

\bibliographystyle{IEEEtran}
  \bibliography{bmc_article}
\end{document}